\documentclass[fmeq,reqno]{amsart}

\usepackage{graphicx}
\usepackage{geometry}
\newtheorem{theorem}{Theorem}[section]

\newtheorem{corollary}[theorem]{Corollary}

\theoremstyle{definition}
\newtheorem{definition}[theorem]{Definition}

\theoremstyle{remark}

\numberwithin{equation}{section}

\newcommand{\N}{\mathbb{N}}

\newcommand{\G}{\mathbb{G}}
\newcommand{\p}{\mathbb{P}}
\newcommand{\M}[1]{\mathbb{#1}}

\newcommand{\mc}[1]{\mathcal{#1}}
\newcommand{\sieb}{\mathcal{S}}

\newcommand{\mood}[1]{(\operatorname{mod} #1)}
\newcommand{\opn}[1]{\operatorname{#1}}

\newcommand{\prodN}[2]{ \prod_{#1 = 3}^#2}

\newcommand{\abs}[1]{\lvert#1\rvert}

\newcommand{\om}[1]{$\omega_{#1}$}

\newcommand{\set}[1]{\lbrace #1\rbrace}

\makeindex

\begin{document}

\title{On the Infinitude of the Twin Primes}

\author{Berndt Gensel}
\address{Prof. Dr. B. Gensel, Spittal on Drau, Austria}
\email{dr.berndt@gensel.at}


\date{\today}

\maketitle

\begin{abstract}
We present a novel approach to the Twin Prime Conjecture, basing on the $6x \pm 1$ representation of primes. By defining so-called twin prime generators $x \in \N$, for which both $6x - 1$ and $6x + 1$ are prime, we reformulate the conjecture into the existence problem of such $x$.

Using admissible residue classes modulo products of small primes and an adapted Selberg sieve, we partition the natural numbers into structured intervals $\mc{A}_n$, where the maximal possible prime divisor of $6x \pm 1$ is fixed. Within each $\mc{A}_n$, we apply the sieve to estimate the number of generator candidates that escape all local obstructions.

Due to the \emph{parity problem} we cannot solve the problem with a Selberg sieve. It requires other sieves or methods. The author is searching for them and invites all interested people to help.
\end{abstract}

\section{Introduction}
The distribution of prime numbers remains one of the central topics in number theory, with particular interest in bounded differences between primes (e.g., Brun’s theorem~\cite{brun}). Since Zhang’s breakthrough~\cite{zhang} on bounded prime gaps, and subsequent improvements by Tao~\cite{tao} and Maynard~\cite{maynard}, it is known that there exist infinitely many pairs of primes with bounded differences. Yet the \emph{Twin Prime Conjecture}, which asserts the infinitude of pairs $(p, p+2)$ with both entries prime, remains open.

In this paper, based on a previous work by the author~\cite{gensel}, we take a structural approach to the problem by examining numbers $x$ for which both $6x - 1$ and $6x + 1$ are prime, so-called \emph{twin prime generators}. Specifically, we construct intervals $\mc{A}_n := [\xi_n, \xi_{n+1})$ with $\xi_n := \frac{p_n^2 - 1}{6}$, such that each $\mc{A}_n$ contains all $x$ whose associated twin pair lies between two consecutive prime squares $p_n^2$ and $p_{n+1}^2$. Here all composite numbers have divisors $\le p_n$.
We define the \om{p_n}-\emph{numbers} as integers $x$, for which $6x\pm 1$ are coprime to all primes $5\le p \le p_n$ and can prove, 
that the \om{p_n}-numbers are representatives of the admissible residue classes $\mathbb{Z} / P_{3,n} \mathbb{Z}$ with respect to two exclusions and
that those in $\mc{A}_n$ are \emph{twin prime generators}.

We apply the Selberg sieve~\cite{selberg} to these intervals, aiming to show that $\mc{A}_n$ contains 
more than one \om{p_n}-numbers. But since the \emph{parity problem} avert the application of sieve methods for lower bounds, we cannot solve the problem with a Selberg sieve. It requires other sieves or methods.

The aim of this paper is to offer the mathematical community the approach to prove the Twin Prime Conjecture via \om{p_n}-numbers.
\bigskip

\section{Twin Prime Generators}

\begin{definition}\label{D-twinprimegenerators}
Let
\begin{align*}
&\N \mbox{ be the set of the positive integers},\\
&\p \mbox{ the set of the prime numbers}, ~ \p^* \mbox{ primes}  \geq 5, \\
&\p_{-}  = \lbrace p \in \p^* \mid p \equiv -1  \mood{6} \rbrace,  
~ \p_{+}  = \lbrace p \in \p^* \mid p \equiv +1  \mood{6} \rbrace  \end{align*}
and
\[
\kappa(n):=
\left\lfloor\frac{n+1}{6}\right\rfloor\mbox{ for } n\in\N
\]
the \textbf{generator} function of the pair $(6\kappa(n)-1,6\kappa(n)+1)$. 
If a pair $(6x-1,6x+1)$ is a twin prime, then we call $x$ a \textbf{twin prime generator} and denote by
\[
\G:=\lbrace x\in\N\mid 6x-1\in\p_-, 6x+1\in\p_+\rbrace
\]
the set of all twin prime generators. Hence the pair $(5,7)$ is the least twin prime in our framework.
\end{definition}

In order to transfer the searching for twin primes to the level of their generators we need a criterion for checking a natural number to be a twin prime generator. 

\begin{theorem}\label{T-twinprimegenerators}
A number $x$ is a \textbf{twin prime generator}, a member of $\G$, if and only if there is \textbf{no} $p \in \p^*$ with $p<6x-1$ that one of the following congruences fulfills:
\begin{align}
x &\equiv - \kappa(p)\mood{p} \label{E-kappa-}
\\
x &\equiv + \kappa(p)\mood{p} \label{E-kappa+}
\end{align}
\end{theorem}
\begin{proof}
At first we assume, that there is a prime $p\in\p^*$ with $p<6x-1$ such that \eqref{E-kappa-} or \eqref{E-kappa+} is valid. There are two cases.
\begin{enumerate}
\item{
$p \in \p_-$ , which means $p=6\kappa(p)-1$:

If \eqref{E-kappa-} is true, then there is an $n \in \M{N}$ with
\begin{align*}
x &= -\kappa(p)+n \cdot (6\kappa(p)-1)\\
6x &= -6\kappa(p)+6n \cdot (6\kappa(p)-1)\\ 
6x+1 &=-6\kappa(p)+6n \cdot (6\kappa(p)-1)+1\\
&=(6n-1)(6\kappa(p)-1)\\
&=(6n-1)\cdot p\\
\implies 6x+1 &\equiv 0 \mood{p} \implies x \notin \G
\end{align*}
For \eqref{E-kappa+} the proof will be done with $6x-1$:
\begin{align*}
6x-1 &=6\kappa(p)+6n \cdot (6\kappa(p)-1)-1\\
&=(6n+1)(6\kappa(p)-1)\\
&=(6n+1)\cdot p\\
\implies 6x-1 &\equiv 0 \mood{p} \implies x \notin \G
\end{align*}
}
\item{
$p \in \p_+$, which means $p=6\kappa(p)+1$:

We proceed similarly with \eqref{E-kappa-} and $6x-1$ as well as \eqref{E-kappa+} and $6x+1$:
\begin{align*}
6x-1 &=(6n-1)(6\kappa(p)+1) \implies 6x-1 \equiv 0 \mood{p}\\
6x+1 &=(6n+1)(6\kappa(p)+1) \implies 6x+1 \equiv 0 \mood{p}
\end{align*}
}
\end{enumerate}
This implies $x \notin \G$, if the congruences \eqref{E-kappa-} or \eqref{E-kappa+} are valid. They cannot both be true, as they are mutually exclusive. \\

If on the other hand $x \notin \G$, then $6x-1$ or $6x+1$ is no prime. 
Let $6x-1 \equiv 0 \mood{p}$ for any $p \in \p_-$. Then we have
\begin{align*}
6x-1 &\equiv 0 \mood{p}\equiv p \mood{p}\\
&\equiv (6\kappa(p)-1) \mood{p}\\
6x &\equiv 6\kappa(p) \mood{p}\\
x &\equiv \kappa(p) \mood{p}\implies \eqref{E-kappa+}.
\end{align*}
For any $p \in \p_+$ we have
\begin{align*}
6x-1 &\equiv -p \mood{p}\\
&\equiv -(6\kappa(p)+1) \mood{p}\\
6x &\equiv -6\kappa(p) \mood{p}\\
x &\equiv -\kappa(p) \mood{p}\implies \eqref{E-kappa-}.
\end{align*}
The other both cases we can handle in the same way. Therefore either \eqref{E-kappa-} or \eqref{E-kappa+} is valid if $x \notin \G$.
\end{proof}

For a composite number $x$ there is at least one factor $\leq\sqrt{x}$. All these we call the \emph{small divisors}.
If we consider, that a small divisor of a composite number $6x+1$ is less or equal to $\sqrt{6x+1}$, then $p$ in the congruences \eqref{E-kappa-} and \eqref{E-kappa+} can be limited by 
\[
\hat{p}(x)= \opn{max}(p \in \p^* \mid p \leq \sqrt{6x+1}).
\]

With $p_n$ as the $n$-th prime number\footnote{It is $p_1=2$}
and $\pi(z)$ as the number of primes $\leq z$ we have with 
\begin{align}\label{E-provable_system}
x \equiv  \kappa(p_n) \mood{p_n}
\mbox{ or }
x \equiv \left( p_n - \kappa(p_n)\right) \mood{p_n}
\end{align}
for $3 \leq n \leq \pi\left(\hat{p}(x)\right)$ a provable system of criteria to exclude all numbers $x \geq 4$ being no twin prime generators.
Since the moduli are primes, the criteria are independent among each other.
\medskip

\section{The \om{p_n}-numbers}

\begin{definition}[Origin]\label{D-origin}
Let
\begin{align}\label{E-origin}
\xi_n := \min(x \in \N \mid \hat{p}(x) = p_n).
\end{align}
It is the first integer $x$, for that the modulus $p_n$ could be a small divisor of $6x\pm 1$. Therefore we denote it as the \textbf{origin} of the modulus $p_n$.
\end{definition}

Let $p_n \leq \hat{p}(x) \leq \sqrt{6x+1}$ and therefore $p_n^2 \leq 6x+1$. Then 
$\dfrac{p_n^2-1}{6}$
is the least number that meets this relation. Comparing with \eqref{E-origin} we get
\begin{align}\label{E-origin2}
\xi_n=\frac{p_n^2-1}{6}.
\end{align}
It is easy to prove, that for every integer $n\geq 3$ holds, that $\xi_n$ is an integer divisible by $4$. 

For every natural number $x$ in the interval
\begin{align}\label{E-A-section}
\mc{A}_n := \lbrace x\in\N\mid \xi_n\leq x<\xi_{n+1}\rbrace 
\end{align}
$\hat{p}(x)$ persists constant on the value $p_n$.
It follows the following Corollary. 

\begin{corollary}\label{C-noDivs>pn}
All numbers $6x\pm 1$ with $x\in\mc{A}_n$ have no small divisors $>p_n$, since $\hat{p}(x)=p_n$ for all $x \in \mc{A}_n$.
\end{corollary}

The length of this interval
will be denoted as $d_n$.
It is depending on the distance between successive primes. 
Since they can only be even, we have with $a = 2,4,6, \ldots$
\begin{align}
d_n & = \frac{(p_n+a)^2-1}{6} - \frac{p_n^2-1}{6}
 = \frac{2ap_n+a^2}{6}\nonumber\\
& = \frac{a}{3}(p_n + \frac{a}{2})
 \geq \frac{2}{3}(p_n+1)>\frac{2p_n}{3}.\label{E-dn1}
\end{align}

Let
\begin{align*}
P_{3,n} &:= \prodN{k}{n}p_k = 5\cdot 7\cdot \ldots \cdot p_n.
\end{align*}

\begin{definition}[\om{p_n}-Numbers]\label{D-omega_numbers}
A positive integer $x$ will be called an \textbf{\om{p_n}-number}, if both $6x-1$ and $6x+1$ are coprime
\footnote{
Then is $(36x^2-1, P_{3,n})=1.$
}
to $P_{3,n}$. Then \eqref{E-provable_system} is for no modulus $p_3,\ldots,p_n$ fulfilled.\\
Vice versa we speak about a \textit{non}-\om{p_n}-\textit{number} $x$, if there is a prime $q \leq p_n$ with $q\mid (36x^2-1)$.
\end{definition}

\begin{theorem}[Main Theorem]\label{T-TPGinAsection} 
Let $x$ be an \om{p_n}-number as a member of $\mc{A}_n$. Then $x$ is a \textbf{twin prime generator}.
\end{theorem} 
\begin{proof}
The claim is valid because by virtue of Definition \ref{D-omega_numbers} 
$6x-1$ as well as $6x+1$ are coprime to $5,7,\ldots,p_n$
and $6x\pm 1$ have by Corollary \ref{C-noDivs>pn} no small divisors that are greater than $p_n$, since $\hat{p}(x) = p_n$ for all members of $\mc{A}_n$.
Hence $6x - 1$ as well as $6x + 1$ must be primes and
$x$ is a twin prime generator.
\end{proof}
\medskip

\section{Admissible Residue Classes}

Consider the residue class ring $\mathbb{Z} / P_{3,n} \mathbb{Z}$. According with Theorem \ref{T-twinprimegenerators} we define the exclusion sets:
\[
\mc{X}_{p_k} := \{ \kappa(p_k), p_k-\kappa(p_k) \} \quad \text{for } k = 3, \ldots, n.
\]

By virtue of the Chinese Remainder Theorem (cf. \cite[Chapter 2, \S~ 2.2]{bundschuh})
the admissible residue classes modulo $P_{3,n}$ with regard to the exclusions are:
\[
\mc{R}_n := \{ r_1, r_2, \ldots, r_m \} = \prod_{k=3}^n \left( \mathbb{Z} / p_k \mathbb{Z} \setminus \mc{X}_{p_k} \right),
\]
with cardinality:
\[
m = \#\mc{R}_n = \prod_{k=3}^n (p_k - 2),
\]
and density:
\begin{align}\label{E-eta}
\eta(p_n) := \frac{\#\mc{R}_n}{P_{3,n}} = \prod_{k=3}^n \left(1 - \frac{2}{p_k}\right).
\end{align}

With Mertens' Theorems can be demonstrated:
\begin{align}\label{E-limEta}
\frac{c_1}{\log^2 p_n}\leq \prod_{k=3}^n \left(1-\frac{2}{p_k}\right)\leq\frac{c_2}{\log^2 p_n}
\text{, implying }\lim_{n\rightarrow\infty} \eta(p_n) = 0.
\end{align}

\begin{theorem}\label{T-omegaARC}
The $\omega_{p_n}$-numbers are representatives of the admissible residue classes modulo $P_{3,n}$ with regard to the exclusions:
\[
\Omega_n := \{ x \in \mathbb{N} \mid x \bmod p_k \notin \mc{X}_{p_k}, k = 3, \ldots, n \} = \mc{R}_n \cap \mc{P}_n,
\]
where $\mc{P}_n := \{ 1, 2, \ldots, P_{3,n} \}$.
\end{theorem}
\begin{proof}
An $x \in \mc{R}_n \cap \mc{P}_n$ satisfies $x \bmod p_k \notin \mc{X}_{p_k}$ for $k = 3, \ldots, n$, which by Theorem \ref{T-twinprimegenerators} implies $(36x^2 - 1, P_{3,n}) = 1$. Thus, $x$ is an $\omega_{p_n}$-number due to Definition \ref{D-omega_numbers}.
\end{proof}
\medskip

\section{Why a Selberg Sieve Fails}\label{S-Sieve}
By virtue of Theorem \ref{T-TPGinAsection}, all \om{p_n}-numbers in $\mc{A}_n$ are twin prime generators, with the cardinality $\abs{\mc{A}_n} = d_n$ (see \eqref{E-dn1}). Therefore we would choose $\mc{A}_n$ by virtue of \eqref{E-A-section} as the sieve set and would define:
\begin{align}
\sieb(\mc{A}_n, \p_n) &:= \abs{\set{x\in\mc{A}_n\mid x \bmod p_k \notin \mc{X}_k \text{, for } k=3,\ldots,n }}\nonumber\\
&= \sum_{\substack{x \in \mc{A}_n \\ (36x^2 - 1, P_{3,n}) = 1}} 1,\label{E-SieveFunction}
\end{align}
where $(36x^2 - 1, P_{3,n})=\opn{gcd}(36x^2 - 1, P_{3,n})$ is the \emph{greatest common divisor}.

Let
\[
\p_n := \set{ p \in \M{P}^* \mid 5 \leq p \leq p_n }
\]
and $\rho(p)$ be the number of solutions to $36x^2 \equiv 1 \mood{p}$. Since $36x^2 - 1 \equiv 0 \mood{p}$ is equivalent to $(6x - 1)(6x + 1) \equiv 0 \mood{p}$, there are exactly two solutions modulo $p$ for all primes $p \geq 5$ (see \cite[Theorem 7.3]{HardyWright}.
Using appropriate Selberg sieve weights $\lambda_d$ for $d \mid P_{3,n}$ (cf. Subsection \ref{S-weights}), we would have:
\begin{align}\label{E-sieveexpression}
\sieb(\mc{A}_n, \p_n) \geq \abs{\mc{A}_n} \cdot V_n - R_n,
\end{align}
where:
\begin{align*}
V_n &:= \prod_{p \in \p_n} \left(1 - \frac{\rho(p)}{p}\right)
\\
&\text{and}\\
R_n &:= \sum_{d_1, d_2 \mid P_{3,n}} \abs{\lambda_{d_1}} \cdot \abs{\lambda_{d_2}} \cdot \abs{\mc{A}_{d_1, d_2}}, \\ 
\mc{A}_{d_1, d_2} &:= \set{ x \in \mc{A}_n \mid 36x^2 \equiv 1 \mood{[d_1, d_2]} }.
\end{align*}
with $[d_1, d_2]=\opn{lcm(d_1,d_2)}$ is the \emph{least common multiple}.
$\abs{\mc{A}_n}\cdot V_n$ is the \emph{main term} and $R_n$ the \emph{error term}.

Since $\rho(p) = 2$, the main term is determined on the density of admissible residue classes (see \eqref{E-eta}):
\begin{align*}
V_n := \prod_{p \in \p_n} \left(1 - \frac{\rho(p)}{p}\right) = \prod_{p \in \p_n} \left(1 - \frac{2}{p}\right) = \eta(p_n).
\end{align*}

But such a sieve would have the \emph{parity problem}, which means, that the sieve cannot differentiate between numbers with even number of prime factor and numbers with odd number of prime factors. 
\medskip

This means, it requires other sieves or methods. The mathematical community is invited to use this shown approach for solving the Twin Prime Conjecture. The author will do it further.


\end{document}